\documentclass[leqno,11pt]{article}

\usepackage{amssymb}
\usepackage{amsthm}
\usepackage{eucal}
\usepackage{amsmath,amsxtra}
\usepackage{graphicx}

\makeatletter

\skip\footins 24pt plus 4pt minus 2pt

\makeatother

\newcommand\bigcheck[1]{#1 \raise1ex\hbox{$\hspace{-1ex}{}^\vee$}}
\newcommand\sucheck[1]{#1 \raise0.5ex\hbox{$\hspace{-1ex}{}^\vee$}}

\newcommand{\checkvarphi}{\smash{\raisebox{-0.9ex}{\kern0.4ex
    \Huge$\check{\smash{\raisebox{0.5ex}{\kern-0.3ex\normalsize
          $\varphi$}}}$}}}

\newcommand{\smcheck}[1]{\smash{\raisebox{-0.8ex}{\kern0.4ex
    \LARGE$\check{\smash{\raisebox{0.5ex}{\kern-0.2ex\normalsize
          $#1$}}}$}}}

\newcommand\subH{{\raise-1ex\hbox{$\scriptstyle H$}}}

\newcommand\barsubH{{\raise-1ex\hbox{$\left| \scriptstyle H \right|$}}}
\newcommand\barsubHo{{\raise-1ex\hbox{$\left| \scriptstyle H_0 \right|$}}}

\makeatletter
\renewcommand\section{\@startsection {section}{1}{\z@}%
                                   {-3.5ex \@plus -1ex \@minus -.2ex}%
                                   {2.3ex \@plus.2ex}%
                                   {\normalfont\normalsize\bfseries}}
\renewcommand\subsection{\@startsection{subsection}{2}{\z@}%
                                     {-3.25ex\@plus -1ex \@minus -.2ex}%
                                     {1.5ex \@plus .2ex}%
                                     {\normalfont\normalsize\bfseries}}
\renewcommand\subsubsection{\@startsection{subsubsection}{3}{\z@}%
                                     {-3.25ex\@plus -1ex \@minus -.2ex}%
                                     {1.5ex \@plus .2ex}%
                                     {\normalfont\normalsize\bfseries}}

\@addtoreset{equation}{section}
\makeatother

\newtheorem*{theorem*}{Support theorem}

\makeatletter
\def\iint{\DOTSI\protect\ints@\tw@}
\def\iiint{\DOTSI\protect\ints@\thr@@}
\def\iiiint{\DOTSI\protect\ints@{4}}
\def\idotsint{\DOTSI\protect\ints@\z@}

\def\intkern@{\mkern-6mu\mathchoice{\mkern-3mu}{}{}{}}
\let\DOTSI\relax
\let\ilimits@\displaylimits

\def\ints@#1{%
  \mkern-7mu\mathchoice{\mkern-2mu}{}{}{}%
  \mathop{\mkern7mu\mathchoice{\mkern2mu}{}{}{}%
    \intop\ifnum#1=\z@\intdots@
    \else\intkern@\fi
    \ifnum#1>\tw@\intop\intkern@\fi
    \ifnum#1>\thr@@\intop\intkern@\fi
    \intop
  }\ilimits@
}
\makeatother

    {
    \begin{list}{}
        { \settowidth{\labelwidth}{}
          \leftmargin=0pt }}
    {\end{list}}

\begin{document}

\begin{center}
\large{\textbf{The Inversion of the X-ray Transform on a Compact
    Symmetric Space}}
\end{center}

\begin{abstract}
The X-ray transform on a compact symmetric space $M$ is here
inverted by means of an explicit inversion formula.  The proof
uses the conjugacy of the minimal closed geodesics in~$M$ and of
the maximally curved totally geodesic spheres in~$M$, proved in
[2b].

\end{abstract}

 The full paper to appear in the Journal of Lie Theory.

\end{document}